%% file: main.tex
\documentclass[11pt]{article}
 \usepackage{geometry}
\usepackage{graphicx}

\usepackage{amsmath, amsfonts, amsthm}
\usepackage{mathtools} 

\usepackage{hyperref}
\usepackage[noabbrev,capitalize]{cleveref} 

\usepackage{graphicx}
\usepackage{subcaption}
\usepackage[font=scriptsize,labelfont=bf,compatibility=false]{caption}
\graphicspath{{./images}}
\usepackage{float}

\usepackage[per-mode = symbol]{siunitx}

\usepackage{bm} 

\newcommand{\R}{\mathbb{R}}

\newcommand{\mr}[1]{\mathrm{#1}}

\usepackage{subfiles} 

\begin{document}

    \title{Efficient design of continuation methods for hyperbolic transport problems in porous media}     
    \author{
        Peter von Schultzendorff \thanks{Department of Mathematics, University of Bergen, Allégaten 44, 5007 Bergen, Norway. Email: \texttt{peter.schultzendorff@uib.no}.}
        \and
        Jakub Wiktor Both \thanks{Department of Mathematics, University of Bergen, Allégaten 44, 5007 Bergen, Norway. Email: \texttt{jakub.both@uib.no}.}
        \and
        Jan Martin Nordbotten \thanks{Department of Mathematics, University of Bergen, Allégaten 44, 5007 Bergen, Norway. NORCE, Nygårdsgaten 112, 5008 Bergen, Norway. Email: \texttt{jan.nordbotten@uib.no}.}
        \and
        Tor Harald Sandve \thanks{NORCE, Nygårdsgaten 112, 5008 Bergen, Norway. Email: \texttt{tosa@norceresearch.no}.}
    }
    \date{}
    \maketitle

    \abstract{
        Full-physics modeling of multiphase flow in porous media, e.g., for carbon storage and groundwater management, requires the nonlinear coupling of various physical processes. Industry standard nonlinear solvers, typically of Newton-type, are not unconditionally convergent and computationally expensive. Homotopy continuation solvers have recently been studied as a robust and versatile alternative. They tackle challenging nonlinear problems by first solving a simple auxiliary problem and then tracing a solution curve towards the more complex target problem. Robustness and efficiency of the method depends on the iterative numerical curve-tracing algorithm as well as on careful design of the auxiliary problem. We assess the traceability of the solution curve for different choices of the auxiliary problem.
        For the Buckley--Leverett equation, modeling two-phase flow in one dimension, we exemplarily compare the previously introduced vanishing-diffusion and linear constitutive laws homotopy continuation, and a new approach based on the entropy solution of the problem. This provides insight toward systematically and robustly designing homotopy continuation methods for solving complex multiphase flow in porous media.
    }

    \section{Introduction}\label{sec:introduction}
        \subfile{introduction}

    \section{Buckley--Leverett equation and its numerical solution}\label{sec:model}

\subfile{model}

    \section{Homotopy continuation method \& traceability metrics}\label{sec:hc}

\subfile{hc}

    \section{Numerical analysis}\label{sec:results}
        \subfile{results}

    \subparagraph{Data availability statement}
        The complete source code is publicly accessible at \url{https://doi.org/10.5281/zenodo.18863360}.

    \subparagraph{Use of generative AI statement}
        The authors used \emph{M365 Copilot} to improve clarity and phrasing of the manuscript. \emph{GitHub Copilot} was used to generate unit tests, write code documentation, and optimize the source code for the traceability measures. The authors have reviewed and revised all AI-generated content and take full responsibility for the final work.

    \subparagraph{Acknowledgements}
        The authors acknowledge funding from the Centre of Sustainable Subsurface Resources (CSSR), grant nr. 331841, supported by the Research Council of Norway, research partners NORCE Norwegian Research Centre and the University of Bergen, and user partners Equinor ASA, Wintershall Dea Norge AS, Sumitomo Corporation, Earth Science Analytics, GCE Ocean Technology, and SLB Scandinavia. JWB acknowledges partial support from the FRIPRO project “Unlocking maximal geological CO2 storage through experimentally validated mathematical modeling of dissolution and convective mixing (TIME4CO2)”, grant nr. 355188, funded by the Research Council of Norway.

\input{references}
\end{document}

%% file: introduction.tex
    Modeling fluid flow in subsurface reservoirs involves strongly nonlinear interactions between fluid phases and components governed by constitutive laws. The resulting nonlinear systems of equations often challenge standard nonlinear solvers, which may fail to converge for large time steps, forcing unnecessary time-step reductions that hinder computational efficiency.

    Hyperbolic transport of fluid phases or components, exemplified by the classical Buckley-Leverett equation for one-dimensional two-phase flow~\cite{Pettersen:GrunnkursReservoarmekanikk1990}, highlights typical difficulties for nonlinear solvers. The characteristics of these hyperbolic, in the presence of capillary forces, weakly parabolic, transport problems are governed by the nonlinear flow functions that determine the fraction of total flow carried by the phase or component of interest. The discretized numerical flow functions introduce inflections and kinks that often cause nonconvergence~\cite{Jenny:UnconditionallyConvergentNonlinear2009,Wang:TrustregionBasedSolver2013,Li:NonlinearAnalysisMultiphase2015}. Moreover, degeneracies in the non-saturated regime cause slowly propagating saturation fronts that further slow nonlinear convergence~\cite{Jiang:DissipationbasedContinuationMethod2018}.

	Several strategies aim to improve nonlinear robustness. Damping methods such as \emph{Appleyard chopping} and trust-region methods constrain saturation updates between Newton steps to prevent crossing steep gradients, inflection points, or kinks in the flow functions~\cite{Jenny:UnconditionallyConvergentNonlinear2009,Wang:TrustregionBasedSolver2013,Li:NonlinearAnalysisMultiphase2015}. Hybrid upwinding methods evaluate viscous, buoyancy, and capillary phase mobilities separately, thereby smoothing the flow functions~\cite{Li:NonlinearAnalysisMultiphase2015,Hamon:ImplicitHybridUpwind2016}. However, except for Appleyard chopping, these methods require explicit fractional flow functions, limiting their use with implicit or black-box constitutive laws, such as hybrid physics‑based, data‑driven models~\cite{Li:AccelerationNVTFlash2019}.

    A versatile alternative is the homotopy continuation (HC) method, which has been applied to solve nonlinear problems in multiphase flow in porous media~\cite{Jiang:DissipationbasedContinuationMethod2018,Salinas:VanishingArtificialDiffusion2020,Jakobsen:HomotopyAnalysisLippmann2020,Jiang:EfficientDissipationbasedNonlinear2023,vonSchultzendorff:AdaptiveHomotopyInPreparation} and computational fluid dynamics~\cite{Brown:MonolithicHomotopyContinuation2016,Brown:DesignEvaluationHomotopies2017}. HC constructs a convex combination of a simpler auxiliary problem and the original system resulting in a solution curve connecting both problems. The HC curve is traced using a predictor-corrector (PC) algorithm. If the solution varies smoothly along the curve, each predictor remains within the convergence region of the subsequent corrector, enabling robust solution of highly nonlinear PDEs that challenge traditional solvers. The overall efficiency and robustness of the HC method depend on the PC scheme and on the auxiliary problem, which governs the geometric traceability of the HC curve.

    This work investigates the HC curve traceability for the Buckley--Leverett equation using three auxiliary problems. We compare curvature, which governs the accuracy of first-order predictors, following the approach of~\cite{Brown:DesignEvaluationHomotopies2017}, and Newton convergence to determine admissible PC step sizes. The auxiliary problems are a vanishing diffusion formulation~\cite{Jiang:DissipationbasedContinuationMethod2018,Salinas:VanishingArtificialDiffusion2020,Jiang:EfficientDissipationbasedNonlinear2023},~a linear constitutive-law approach~\cite{vonSchultzendorff:AdaptiveHomotopyInPreparation}, and a new formulation based on the entropy solution of the hyperbolic transport problem. To our knowledge, curve traceability for continuation methods in porous media flow has not been systematically assessed.

    The rest of the article is organized as follows. \cref{sec:model} describes the Buckley–Leverett problem and discretization. \cref{sec:hc} introduces the HC method, auxiliary problems, and traceability metrics. \cref{sec:results} presents the numerical assessment of HC curves.

%% file: model.tex
    For this study, we consider the Buckley--Leverett equation, which describes incompressible and immiscible flow of two fluid phases in a one-dimensional homogeneous porous medium. The one-dimensional setting enables direct visualization and study of the HC curves. Despite its simplicity, the Buckley--Leverett equation captures key complexities that also appear in the higher-dimensional full-physics model, for example, when one fluid displaces another through narrow channels in heterogeneous formations, where the behavior is effectively one-dimensional. Given a domain \(\Omega = [0,L]\) and a time interval \(I = [0,T]\), the mathematical model reads
    \begin{subequations}\label{eq:model:Buckley--Leverett problem}
        \begin{align}
            \frac{\partial S}{\partial t} + \frac{v}{\phi} \frac{\partial f(S)}{\partial x} &= 0 \quad \text{on } \Omega \times (0,T],\\
            S(x, 0) &= S^0, \quad S(0, t) = S_i,
        \end{align}
    \end{subequations}
    Here, the unknown \(S(x,t) \in [0,1]\) describes the wetting phase saturation, the flow function \(f(S)\) represents the ratio of wetting volumetric flow rate to the total flow rate, and \(v\) and \(\phi\) are total volumetric flow and porosity, respectively. 
    By choosing a suitable rescaling of the time \(t\), we can assume that \(\frac{v}{\phi} = 1\). For the common Corey relative permeability model \(k_{r,\mr{w}} = S^{n_\mr{w}}\) and \(k_{r,\mr{n}} = (1 - S)^{n_\mr{n}}\) and a constant viscosity ratio \(M = \frac{\mu_\mr{w}}{\mu_\mr{n}}\), the fractional flow function is given by 
    \begin{equation}
        f(S) = \frac{S^{n_w}}{S^{n_w} + M (1 - S)^{n_n}}
    \end{equation}

    \subsection{Fully implicit finite volume discretization}
        To solve~\eqref{eq:model:Buckley--Leverett problem} numerically, we consider a partitioning of the spatial domain into \(N_K\) control volumes \(\left(x_{k-\frac{1}{2}}, x_{k+\frac{1}{2}}\right)\) of uniform size \(\Delta x\), \(k=1, ..., N_K\), and a partitioning of the temporal domain into time intervals \([t^{n-1}, t^n)\) of size \(\tau^n = t^n - t^{n-1}\), \(n=1, ..., N_n\). We denote by \(S^n_k\) the approximate cell average of the saturation in cell \(K_k\) at time \(t^n\). Discretizing the Buckley--Leverett equation with the implicit Euler method, cell-centered finite volumes (CCFVM), and upstream mobility weighting, we obtain for each discrete time and each control volume the discrete problem
        \begin{equation}\label{eq:model:discretized problem}
            \frac{S_k^{n} - S_k^{n-1}}{\tau^n} + \frac{f(S_{k-1}^{n}) - f(S_k^{n})}{\Delta x} = 0, \quad k = 1, \dots, N_K, n=1, \dots, N_n.
        \end{equation}

    \subsection{Newton's method \& nonlinear solver challenges}\label{subsection:model:newton}
        The nonlinear system of coupled equations at discrete time \(t^n\), denoted by \(\mathbf{F}^n(\mathbf{X}^n) = 0\), is typically solved with Newton's method. Newton's method is also the nonlinear solver of our choice for the corrector step of the HC method, as described in~\cref{sec:hc}. However, we emphasize that the presented HC method and the framework for numerical analysis can be combined with other iterative nonlinear solvers.

        In the purely viscous case, Newton's method is unconditionally convergent for linear, convex, or concave flux functions. However, for S-shaped flux functions Newton may be non-convergent whenever a cell saturation update crosses an inflection point where \(f''(S) = 0\)~\cite{Jenny:UnconditionallyConvergentNonlinear2009}. In the degenerate regime, where one phase invades a medium initially saturated with the other and the initial saturation is close to an immobile state \(S^0 = 0\) or \(S^0 = 1\), nonlinear relative permeabilities cause the saturation front to advance only by one grid cell per Newton step. This severely limits the overall convergence rate~\cite{Jiang:DissipationbasedContinuationMethod2018}.

%% file: hc.tex
    The homotopy continuation (HC) method is a globalization technique that aims at improving robustness and overall convergence speed compared to Newton's method. It has proven effective for solving nonlinear PDEs, in CFD~\cite{Brown:MonolithicHomotopyContinuation2016,Brown:DesignEvaluationHomotopies2017} and porous media flow~\cite{Jiang:DissipationbasedContinuationMethod2018,Salinas:VanishingArtificialDiffusion2020,Jiang:EfficientDissipationbasedNonlinear2023}. Rather than solving the nonlinear system
    \begin{equation}\label{eq:target problem}
        \mathbf{F}: \R^n \to \R^n, \quad \mathbf{F}(\mathbf{X}) = 0
    \end{equation}
    directly, the problem is embedded into the higher-dimensional space \(\R^{n+1}\) using the HC framework. The key idea is to iteratively deform a simpler auxiliary problem 
    \begin{equation}\label{eq:initial problem}
        \mathbf{G}: \R^n \to \R^n, \quad \mathbf{G}(\mathbf{X}) = 0,
    \end{equation}
    into the target problem. The auxiliary problem should be chosen such that its solution can be efficiently obtained with standard nonlinear solvers. The HC method then considers a convex combination with the target problem
    \begin{equation}\label{eq:HC problem}
        \mathbf{H}: \R^{n+1} \to \R^n, \quad \mathbf{H}(\mathbf{X}; \lambda) \coloneqq \lambda \mathbf{G}(\mathbf{X}) + (1 - \lambda) \mathbf{F}(\mathbf{X}) = 0, \quad \lambda \in [0, 1].
    \end{equation}
    For \(\lambda = 1\), the HC problem reduces to the auxiliary problem \(\mathbf{H}(\mathbf{X};1) = \mathbf{G}(\mathbf{X})\), while for \(\lambda = 0\) the target problem \(\mathbf{H}(\mathbf{X};0) = \mathbf{F}(\mathbf{X})\) is recovered. The HC method proceeds by tracing the solution curve
    \begin{equation}
        \mathcal{H}:= \left\{\left(\mathbf{X};\lambda\right) \mid \exists \lambda \in [0,1] \text{ s.t. } \mathbf{H}\left(\mathbf{X};\lambda\right) = 0\right\}
    \end{equation}
    from \(\lambda = 1\) to \(\lambda = 0\). For the numerical solution of PDEs, a predictor-corrector (PC) algorithm is commonly used for curve tracing~\cite{Brown:DesignEvaluationHomotopies2017,Jiang:DissipationbasedContinuationMethod2018}, while alternative algorithms exist~\cite{Allgower:IntroductionNumericalContinuation2003}. The PC algorithm discretizes \(\lambda\) as \(\lambda^0 = 1,\lambda^1,\lambda^2,\dots\to 0\). At each HC step \(i\), the predictor provides an initial guess for solving the subproblem \(\mathbf{H}\left(\mathbf{X},\lambda^i\right) = 0\). The corrector refines the initial guess by solving the subproblem up to the desired precision. A common choice is to use a first-order Euler predictor step and Newton's method as the corrector~\cite{Allgower:IntroductionNumericalContinuation2003,Brown:DesignEvaluationHomotopies2017}.

    \subsection{HC design criteria \& auxiliary problem}
        Robustness and efficiency of the HC method hinge on how accurately the PC algorithm traces the curve. Robustness requires corrector convergence at each continuation step, necessitating accurate predictor steps and/or sufficiently short step sizes \(\Delta \lambda = \lambda^{i+1} - \lambda^i\). Efficiency demands minimizing the total number of PC steps. These challenges can be addressed in two complementary ways. On the algorithmic level, the PC method can be improved through higher-order predictors, better nonlinear solvers in the corrector, or adaptive stopping criteria to reduce overall computational cost. Conversely, geometric traceability of the HC curve can be improved by optimizing the auxiliary problem.

        Following the design criteria for HC outlined in~\cite{Brown:DesignEvaluationHomotopies2017}, the auxiliary problem should be chosen such that:
        \begin{enumerate}
            \setlength{\itemsep}{0pt}
            \item The auxiliary problem \(\mathbf{G}(\mathbf{X}) = 0\) has a unique solution that can be found robustly, ideally unconditionally, using a standard nonlinear solver.
            \item The HC curve can be traced reliably and efficiently using a PC method. This is governed by two key characteristics of the curve:
                \begin{enumerate}
                    \setlength{\itemsep}{0pt}
                    \item The curvature \(\kappa\) of the curve is small, enabling accurate initial guesses to the subproblems \(\mathbf{H}(\mathbf{X};\lambda) = 0\) with low-order predictors.
                    \item The nonlinear subproblems admit large convergence regions, allowing corrector convergence from relatively poor initial guesses. 
                \end{enumerate} 
            \item The HC curve is free of bifurcation points that necessitate small steps along the curve~\cite{Allgower:IntroductionNumericalContinuation2003}. 
            \item The linearized subproblems along the HC curve are computationally tractable and ideally simpler to solve than those of the target problem.
        \end{enumerate}    

        To improve robustness and efficiency of HC for the Buckley--Leverett problem, we focus on the choice of the auxiliary problem with regard to the first two criteria listed above. We compare three auxiliary problems presented in more detail in the following subsections.

        \subsubsection{Vanishing artificial diffusion}
            Suggested for use in porous media flow in~\cite{Jiang:DissipationbasedContinuationMethod2018} and extended in ~\cite{Salinas:VanishingArtificialDiffusion2020,Jiang:EfficientDissipationbasedNonlinear2023}, the vanishing diffusion HC augments the governing equations with an artificial diffusion term
            \begin{subequations}
            \begin{equation}
                \mathbf{H}_\mr{d}(\mathbf{X}) = \mathbf{F}(\mathbf{X}) + \lambda \mathbf{D}(\mathbf{X}).
            \end{equation}
            In the discretized Buckley--Leverett problem, the artificial diffusion term reads
            \begin{equation}
                \beta \frac{\tau^n}{\Delta x} \left(S^{n+1}_{i-1} - 2S^{n+1}_{i} + S^{n+1}_{i+1}\right),
            \end{equation}
            \end{subequations}
            where \(\beta\) is a tunable parameter. The diffusion regularizes the hyperbolic transport problem, by making the auxiliary problem more parabolic. Consequently, at \(\lambda  = 1\) mass can disperse globally across the domain, towards the approximate final destination. As \(\lambda \to 0\), the diffusion term vanishes and the target problem is recovered, but mass only has to move locally. This enables robust convergence if the Newton steps in the corrector do not have to cross inflection or transition points in the flow function. The vanishing diffusion HC has also been shown to speed up convergence in degenerate cases where the initial saturation is close to the immobile saturation \(S^0 = 0\) or \(S^0 = 1\) ~\cite{Jiang:DissipationbasedContinuationMethod2018}. 

            A drawback of the vanishing diffusion approach is that the HC curve formally extends to the fully diffusive limit as \(\lambda \to \infty\), so a finite starting point of the curve must be imposed by prescribing the diffusion strength at \(\lambda = 1\). For the one-dimensional transport problem, a CFL-type argument and numerical tests suggest \(\beta = \omega \max |f'(S)|\), with \(\omega = 2 \cdot 10^{-3}\), whereas for two-dimensional coupled flow and transport a much smaller value \(\omega = 10^{-5}\) is suggested~\cite{Jiang:DissipationbasedContinuationMethod2018}. This highlights the inherent limitation of the artificial diffusion HC that necessitates numerical determination of the initial diffusion strength.
    
        \subsubsection{Linear relative permeabilities}
            Linear relative permeabilities yield a linear, convex, or concave flow function, cf.~\cref{fig:results:case_1:flow function,fig:results:case_2a:flow function}, for which the discretized problem~\eqref{eq:model:discretized problem} is unconditionally convergent~\cite{Jenny:UnconditionallyConvergentNonlinear2009}. In the degenerate regime, linear relative permeabilities also allow the saturation front to advance more than one grid cell per Newton step, improving convergence. This motivates starting the HC with linear relative permeabilities:
            \begin{subequations}
                \begin{align}
                    \mathbf{G_\mr{lin}^n(\mathbf{X}^n)}_k &= \frac{S_k^{,n+1} - S_k^{,n}}{\tau^n} + \frac{f_\mr{lin}(S_{k-1}^{n+1}) - f_\mr{lin}(S_k^{n+1})}{\Delta x}, \\
                    \mathbf{H}_\mr{lin} &= \lambda \mathbf{G}_\mr{lin}(\mathbf{X}) + (1 - \lambda) \mathbf{F}(\mathbf{X}).
                \end{align}
            \end{subequations}
            This HC has previously been used in~\cite{vonSchultzendorff:AdaptiveHomotopyInPreparation}, where it was combined with adaptive stopping criteria.

        \subsubsection{Convex/concave hull}
            For hyperbolic conservation problems, the analytical solution of the Riemann problem depends only on the upper concave or lower convex hull of the flow function, depending on wave direction. Replacing the flow function \(f(S)\) with its convex hull \(f_\mr{c}(S)\) therefore yields a modified Buckley--Leverett problem with the same entropy solution. This motivates an HC design in which the auxiliary problem uses the convex/concave hull. By construction, \(f_\mr{c}(S)\) is globally convex/concave and linear near the immobile saturation of the invading phase. Consequently, the auxiliary problem offers the same benefits as linear relative permeabilities: unconditional convergence for all initial conditions and fast convergence in degenerate regimes.
            \begin{subequations}
                \begin{align}
                    \mathbf{G_\mr{c}^n(\mathbf{X}^n)}_k &= \frac{S_k^{,n+1} - S_k^{,n}}{\tau^n} + \frac{f_\mr{c}(S_{k-1}^{n+1}) - f_\mr{c}(S_k^{n+1})}{\Delta x}, \\
                    \mathbf{H}_\mr{c} &= \lambda \mathbf{G}_\mr{c}(\mathbf{X}) + (1 - \lambda) \mathbf{F}(\mathbf{X})
                \end{align}
            \end{subequations}
            This idea is conceptually similar to the corrected operator splitting for convection-diffusion equations. The convective splitting step places the shock at the correct location, analogous to solving the convex/concave hull problem, while the diffusive splitting step restores the correct shock profile, analogous to solving the target problem~\cite{Karlsen:OperatorSplittingMethods2001}.

    \subsection{HC curve traceability}
        To study curvature and convergence properties of the Newton corrector step along the HC curve we introduce two quantitative metrics. We emphasize that curvature and Newton convergence along the HC curve are highly problem-dependent. They are influenced by the physical model, the geometry and dimensionality, and the discretization of the problem of interest. Therefore, extrapolating results on traceability of HC curves from the Buckley--Leverett problem to higher-dimensional and full-physics simulations is possible only on an abstract conceptual level.

        \subsubsection{HC curvature}
            Given an arclength parametrization \(\mathbf{q}(s) \coloneqq \left(\mathbf{X}(s),\lambda(s)\right)\) of the HC curve \(\mathcal{H}\) satisfying \(\dot{\mathbf{q}}(s) \cdot \dot{\mathbf{q}}(s) = 1\) and assuming the HC curve is traced with steps of constant arclength \(\Delta s\), a Taylor expansion of \(\mathbf{H} \left(\mathbf{X}(s);\lambda(s)\right)\) shows that the prediction error of a first-order Euler predictor scales as \(e \sim s_\mr{tot}^2 \kappa\), where
            \begin{equation}
                \kappa(s) = \sqrt{\ddot{\mathbf{q}}(s) \cdot \ddot{\mathbf{q}}(s)}
            \end{equation}
            is the total curvature of the HC curve and \(s_\mr{tot}\) is the total arclength of the curve~\cite{Brown:DesignEvaluationHomotopies2017}. Although the magnitude of \(s_\mr{tot}^2 \kappa\) depends on the specific problem, comparatively large values indicate that large predictor steps taken along the tangent may yield poor initial guesses for the Newton corrector loop.

        \subsubsection{Newton convergence along the HC curve}
            Determining the region of Newton convergence for complex nonlinear problems is inherently challenging. To the best of the authors’ knowledge, no theoretical result on Newton convergence exists for the discretized one-dimensional hyperbolic transport problem when the flow function is neither linear, convex, nor concave.

            In practice, PC algorithms often use zero-order or first-order predictors~\cite{Allgower:IntroductionNumericalContinuation2003,Brown:DesignEvaluationHomotopies2017}. In both cases, given a previous solution along the curve, \(\mathbf{q}(s)\), the initial guess for the Newton corrector step is taken along the line \(\mathbf{q}_\mr{pred}(s)(\gamma) = \{\mathbf{q}(s) + \gamma \dot{\mathbf{q}}(s) \mid \gamma \in \R^+\}\). To quantify the longest admissible predictor step size along the line, we define
            \begin{subequations}
                \begin{align}
                    r(s) &= \max_r \left\{ r \;\middle|\; 
                    \begin{aligned}
                        &\text{Newton converges for } \mathbf{H}\left(\mathbf{X}; \lambda_\mathrm{pred}(s)(\gamma)\right) \\
                        &\text{from } \mathbf{q}_\mathrm{pred}(s)(\gamma) \, \forall \, 0 < \gamma \leq r
                    \end{aligned}
                    \right\}, \\
                    \intertext{and the normalized measure}
                    \tilde{r}(s) &= \frac{r(s) \dot{\lambda}(s)}{\lambda(s)} \in [0,1].
                \end{align}
            \end{subequations}
            A value of \(\tilde{r}(s) = 1\) indicates that the target system can be solved with an initial guess on the tangent predictor.

            Alternatively, nonlinear complexity along the curve can be characterized by auxiliary dimensionless numbers~\cite{Brown:DesignEvaluationHomotopies2017}. In porous media flow these are the Peclet number, representing the ratio of viscous to capillary forces, the gravity number, representing the ratio of buoyancy to viscous forces~\cite{Li:NonlinearAnalysisMultiphase2015}, and the CFL number, corresponding to the maximum domain of wave propagation~\cite{Jiang:DissipationbasedContinuationMethod2018}.

%% file: results.tex
    We compare curvature and Newton convergence of the linear relative permeability HC, the convex/concave hull HC, and the vanishing diffusion HC with three initial diffusion strengths \(\omega \in \left\{10^{-5}, 2 \times 10^{-3},10^{-1}\right\}\), on four Riemann problems. Constant initial conditions \(S^0 \mid_{[0,0.5]} = S^0_\mr{right}, S^0 \mid_{[0.5,1]} = S^0_\mr{left}\) are prescribed on each half of the domain. Initial conditions, saturation at the inlet \(S_i\) and the viscosity ratio \(M\) are varied as reported in~\cref{tab:parameters}, while the Corey exponents are kept constant \(n_\mr{w} = n_\mr{n} = 2\). The domain is discretized into \(N_K = 100\) cells and a single time step of size \(\tau^1 = 25.0\) is solved. Appleyard damping and physical saturation bounds \(S \in [0,1]\) are enforced after each Newton update in the corrector loop.

    \begin{table}[t]
        \centering
        \caption{Parameters for the Riemann problem test cases (\(n_\text{w} = n_\text{n} = 2\), \(N_K = 100\), \(\tau^1 = 25.0\)).}
        \label{tab:parameters}
        \begin{minipage}{0.45\linewidth}
            \centering
            \begin{tabular}{l p{0.8cm} p{1.2cm} p{0.5cm}}
            \hline\noalign{\smallskip}
            Case & $M$ & $S^0_{L,R}$ & $S_i$ \\
            \noalign{\smallskip}\hline\noalign{\smallskip}
            1   & 1.0  & 0, 1 & 0.0 \\
            2a  & 10.0 & 0, 1 & 0.0 \\
            \noalign{\smallskip}\hline
            \end{tabular}
        \end{minipage}
        \begin{minipage}{0.45\linewidth}
            \centering
            \begin{tabular}{l p{0.8cm} p{1.2cm} p{0.5cm}}
            \hline\noalign{\smallskip}
            Case & $M$ & $S^0_{L,R}$ & $S_i$ \\
            \noalign{\smallskip}\hline\noalign{\smallskip}
            2b  & 10.0 & 1, 0 & 1.0 \\
            2c  & 10.0 & 0.8, 0.2 & 0.8 \\
            \noalign{\smallskip}\hline
            \end{tabular}
        \end{minipage}   
    \end{table}

    Across all test cases, the linear relative permeability HC shows relatively high curvature at \(\lambda = 1\), but curvature consistently decreases along the curve, indicating good overall curve traceability regardless of problem parameters.

\begin{figure}[htbp]
        \centering
        \subcaptionbox{Curvature (Case 1)\label{fig:results:case_1:curvature}}[0.48\linewidth]{\includegraphics[width=\linewidth]{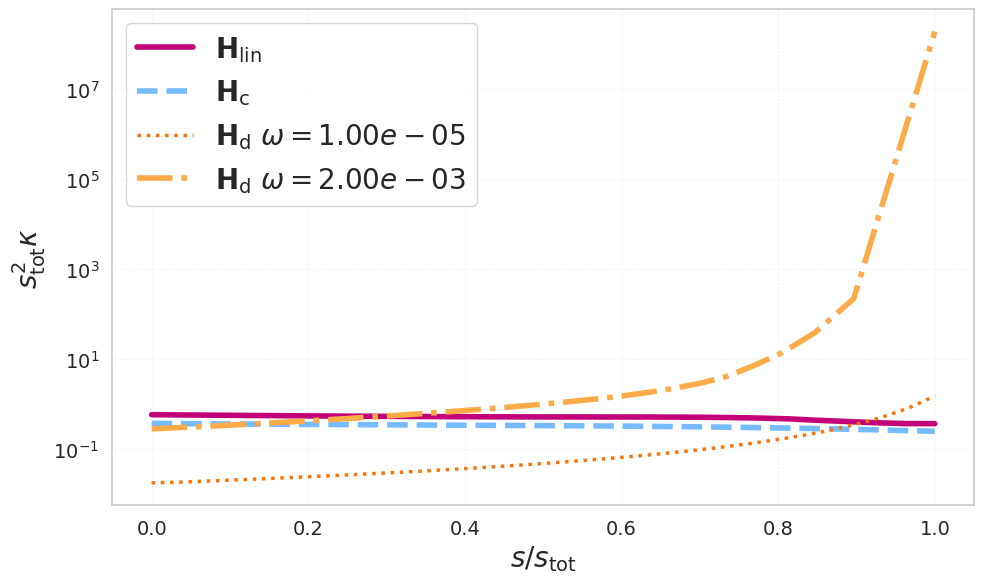}}\hfill
        \subcaptionbox{Metric (Case 1)\label{fig:results:case_1:convergence}}[0.48\linewidth]{\includegraphics[width=\linewidth]{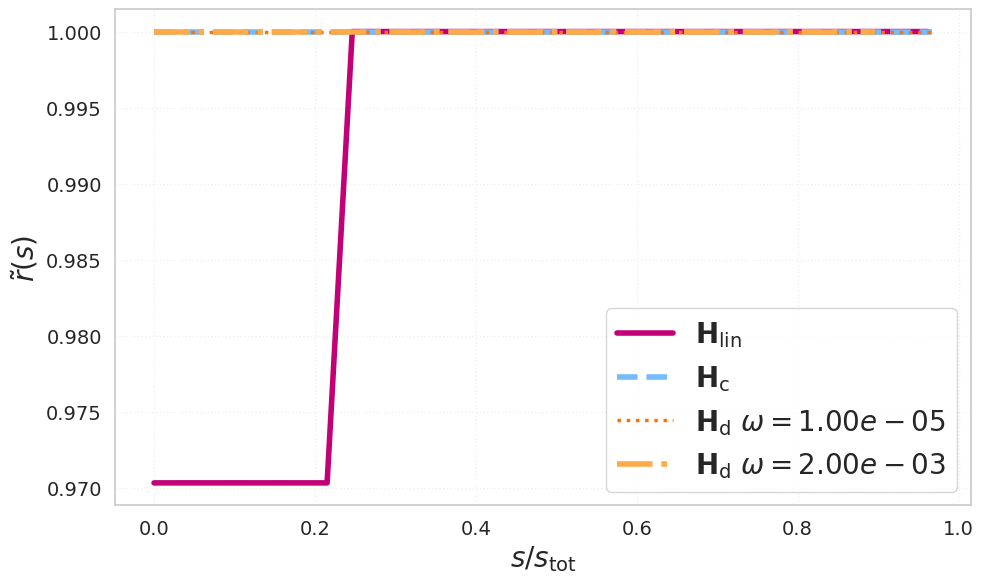}}
        
        \vspace{0.5em}
        
        \subcaptionbox{Curvature (Case 2a)\label{fig:results:case_2a:curvature}}[0.48\linewidth]{\includegraphics[width=\linewidth]{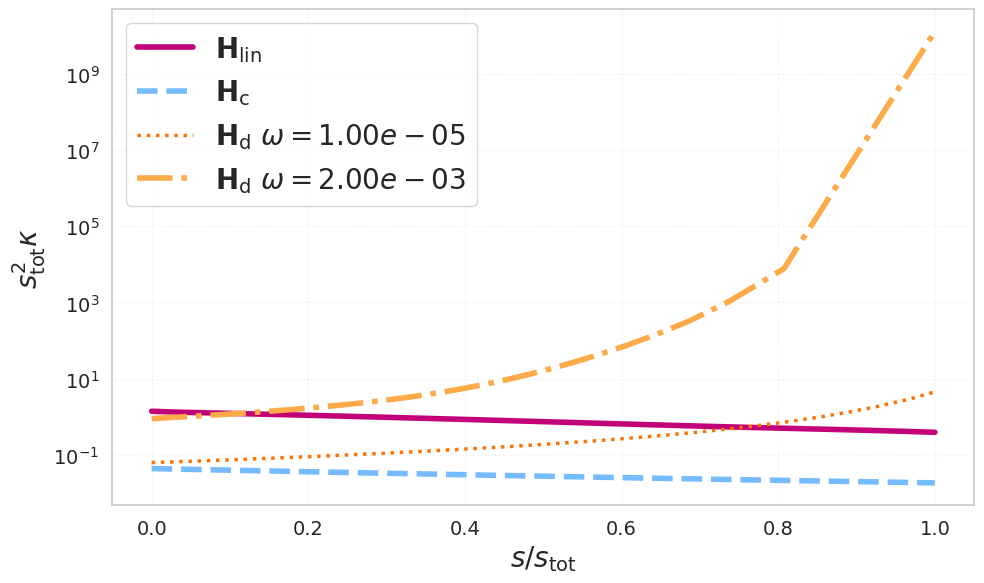}}\hfill
        \subcaptionbox{Metric (Case 2a)\label{fig:results:case_2a:convergence}}[0.48\linewidth]{\includegraphics[width=\linewidth]{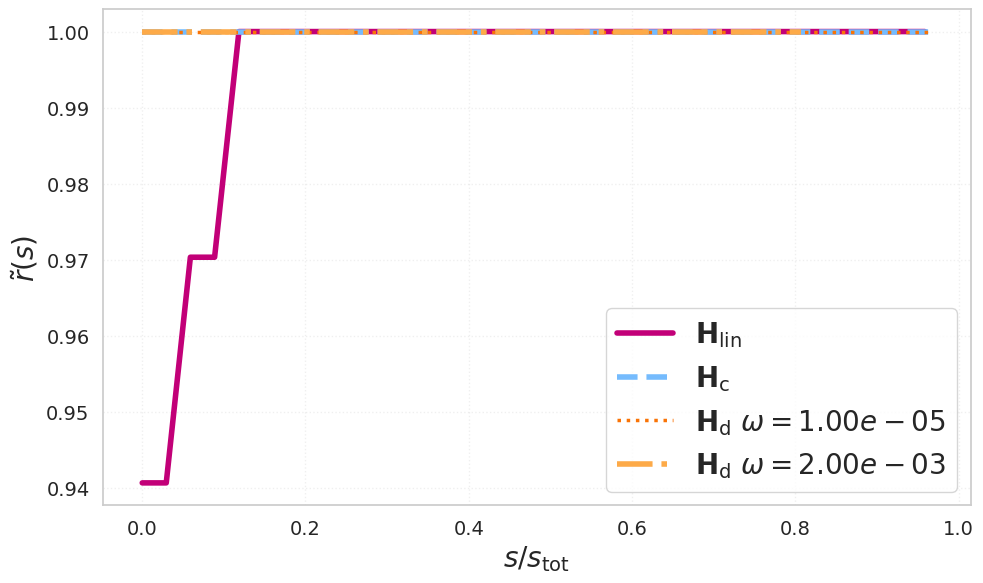}}
        
        \caption{Curvature and relative convergence metrics for Case 1 (\(M=1\)) and Case 2a (\(M=10\)).}
        \label{fig:results:combined_1_2a}
    \end{figure}

\begin{figure}[htbp]
        \centering
        \includegraphics[width=0.40\linewidth]{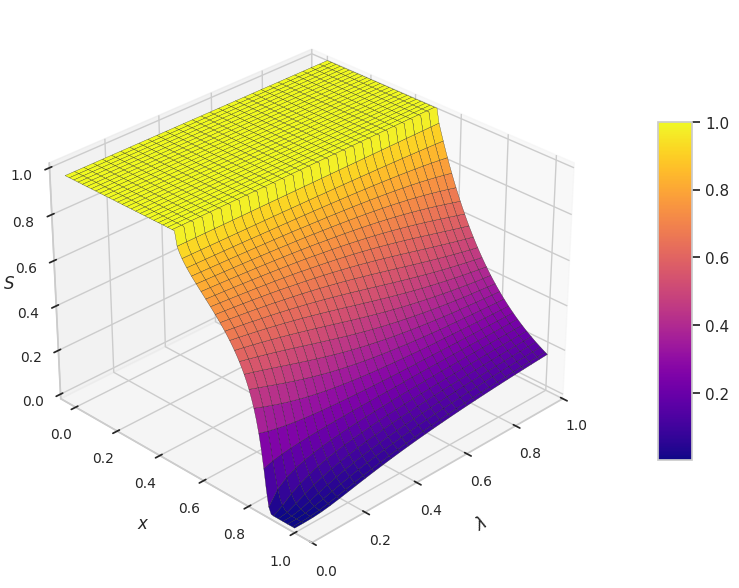}
        \caption{Solution curve for Case 2b using the convex/concave hull HC.}
        \label{fig:results:case_2b:solution curve}
    \end{figure}

    When the convex/concave hull closely matches the target flow function, it achieves the lowest curvature among all HC designs, cf.~\cref{fig:results:combined_1_2a}. However, when the invading phase has the larger viscosity, the linear part of the convex/concave hull deviates significantly from the target flow function, increasing curvature. Inspection of the solution curve suggests that this is caused by added numerical diffusion that prevents accurate numerical resolution of the shock at \(\lambda = 1\), cf.~\cref{fig:results:case_2b:solution curve}.

\begin{figure}[htbp]
        \centering
        \subcaptionbox{Curvature (Case 2b)\label{fig:results:case_2b:curvature}}[0.48\linewidth]{
            \includegraphics[width=\linewidth]{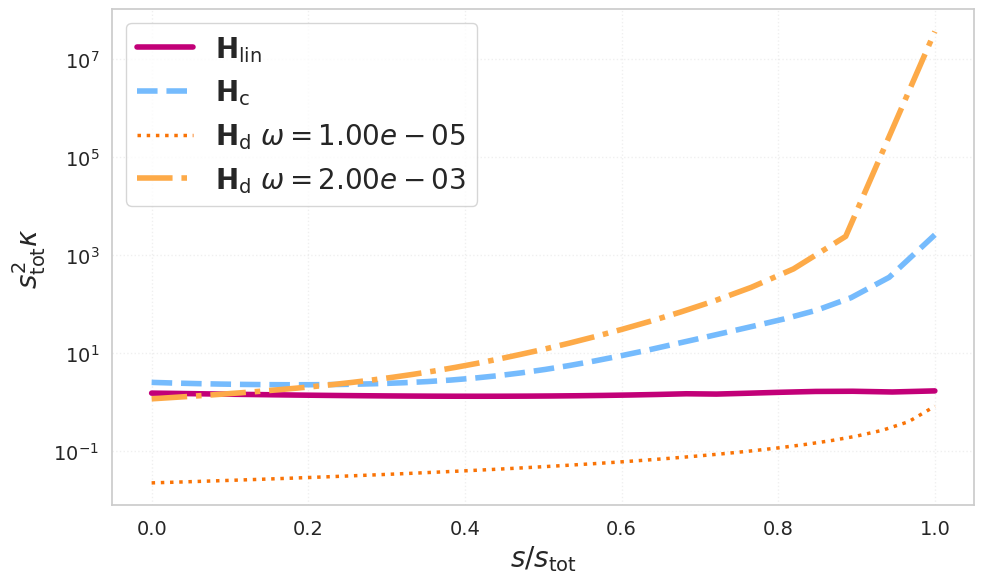}}\hfill
        \subcaptionbox{Metric (Case 2b)\label{fig:results:case_2b:convergence}}[0.48\linewidth]{
            \includegraphics[width=\linewidth]{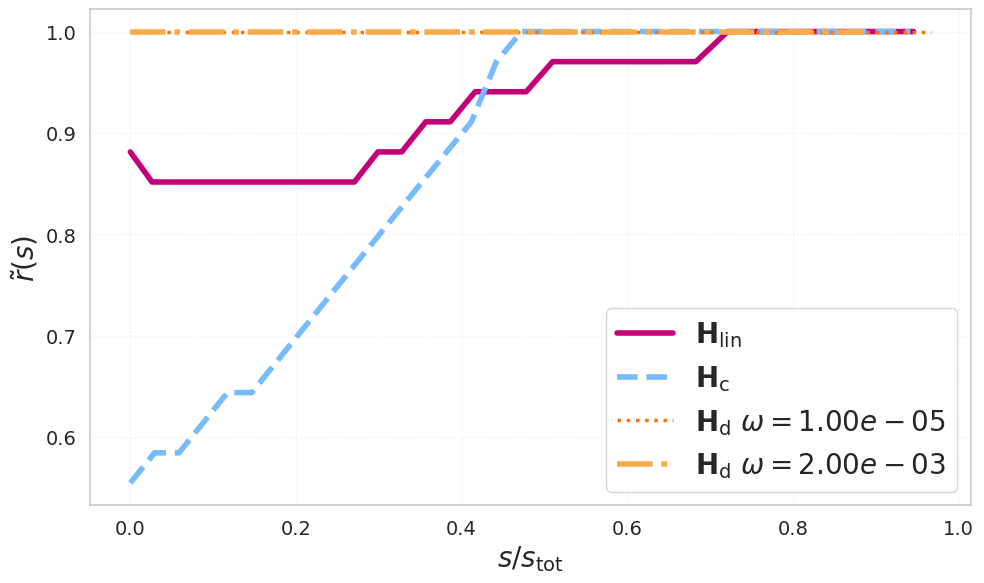}}
        
        \vspace{0.5em} 
        
        \subcaptionbox{Curvature (Case 2c)\label{fig:results:case_2c:curvature}}[0.48\linewidth]{
            \includegraphics[width=\linewidth]{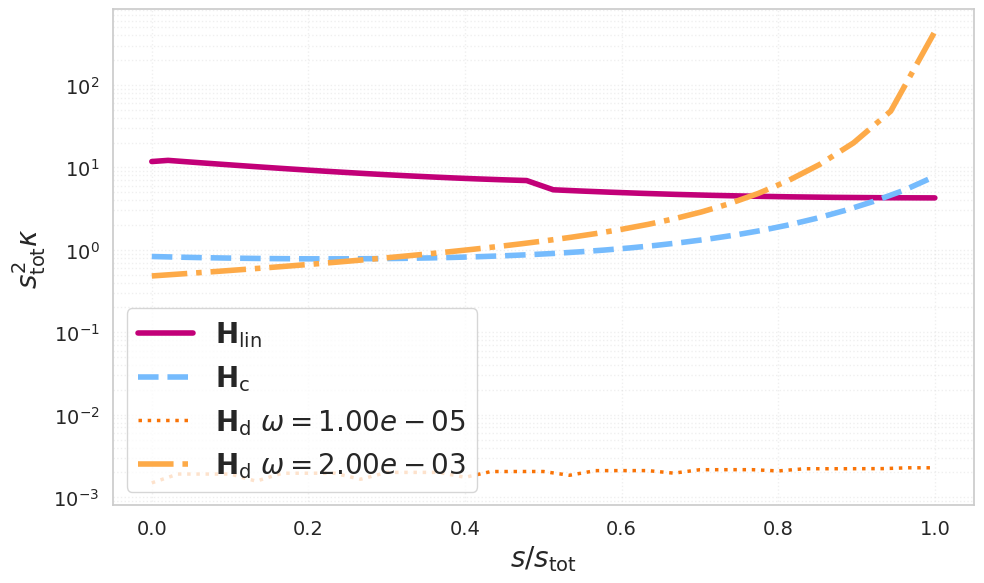}}\hfill
        \subcaptionbox{Metric (Case 2c)\label{fig:results:case_2c:convergence}}[0.48\linewidth]{
            \includegraphics[width=\linewidth]{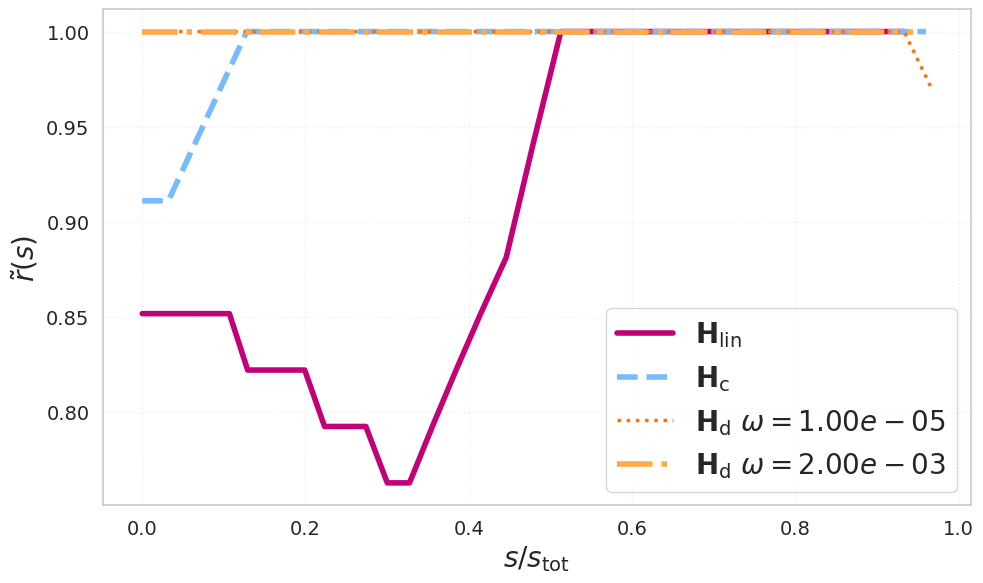}}
        
        \caption{Curvature and relative convergence metrics for Case 2b (\(M=10, S_i=1.0\)) and Case 2c (\(M=10, S_i=0.8\)).}
        \label{fig:results:combined_2b_2c}
    \end{figure}

    The vanishing diffusion HC with \(\omega = 2 \times 10^{-3}\) shows a rapid increase in curvature in the last third of the curve for all problems, cf.~\cref{fig:results:case_1:curvature,fig:results:case_2a:curvature,fig:results:case_2b:curvature,fig:results:case_2c:curvature}. Nevertheless, early-stage curvature remains low enough for \(\tilde{r}(s)\) to attain its maximum along the entire curve. These findings should be interpreted cautiously. Even a slight initial curvature increase reduces \(\tilde{r}(s)\), as seen for the linear relative permeability HC in~\cref{fig:results:case_2b:convergence}. In the non-degenerate case~\cref{fig:results:case_2c:curvature}, the curvature is significantly lower than in the degenerate cases with immobile invading phase saturation~\cref{fig:results:case_1:curvature,fig:results:case_2a:curvature,fig:results:case_2b:curvature}. 

    For all problems, the vanishing diffusion HC with \(\omega = 10^{-1}\) fails to converge at \(\lambda = 1\). Conversely, the HC with \(\omega = 10^{-5}\) converges robustly and achieves the lowest curvature in three of the four cases, indicating easy traceability~\cref{fig:results:case_1:curvature,fig:results:case_2a:curvature,fig:results:case_2b:curvature,fig:results:case_2c:curvature}. However, it is presumably a poor choice for the auxiliary problem, as it is too close to the target problem, as discussed below. The results confirm that \(\omega = 2 \times 10^{-3}\) is indeed a suitable choice for the Buckley--Leverett equation.

    \begin{figure}[htbp]
        \centering
        \subcaptionbox{Case 1\label{fig:results:case_1:flow function}}[0.31\linewidth]{
            \includegraphics[width=\linewidth]{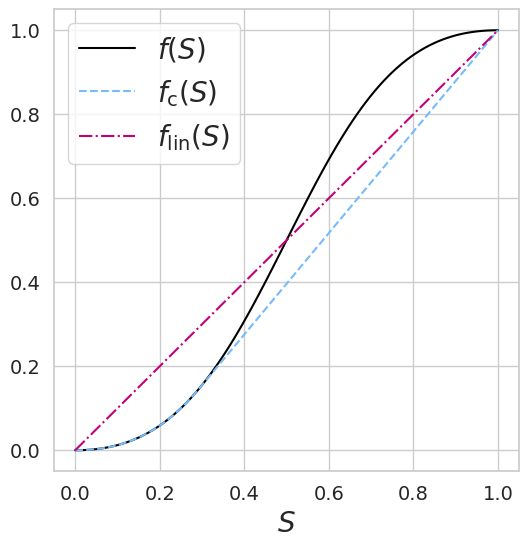}
        }
        \hfill
        \subcaptionbox{Case 2a\label{fig:results:case_2a:flow function}}[0.31\linewidth]{
            \includegraphics[width=\linewidth]{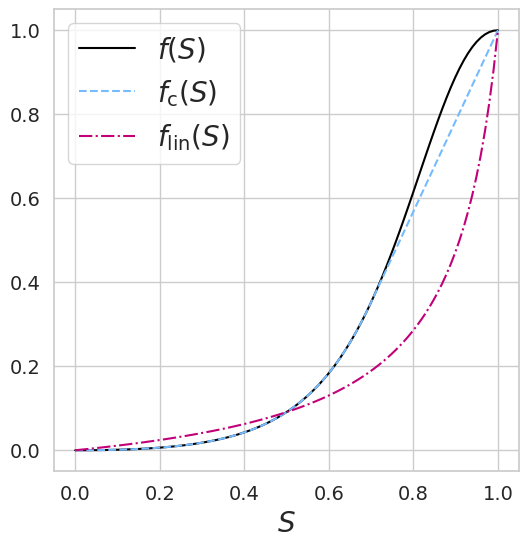}
        }
        \hfill
        \subcaptionbox{Case 2b and 2c\label{fig:results:case_2b:flow function}}[0.31\linewidth]{
            \includegraphics[width=\linewidth]{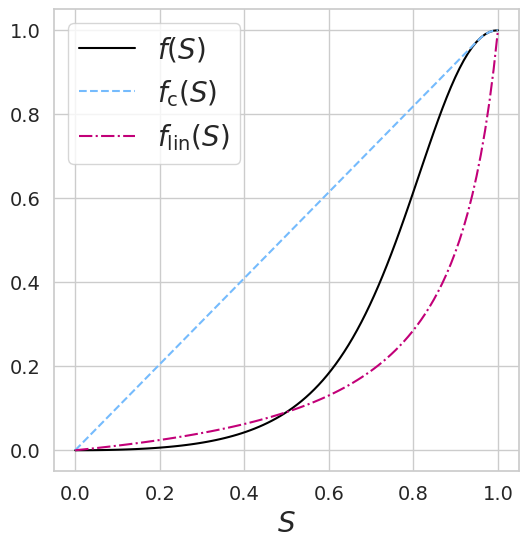}
        }
        \caption{Target, linear relative permeabilities, and convex hull flow functions for the different cases.}
        \label{fig:results:flow functions}
    \end{figure}

\begin{figure}[htbp]
        \subcaptionbox{Linear relative permeabilities HC.}[0.40\linewidth]
            {
                \includegraphics[width=\linewidth]{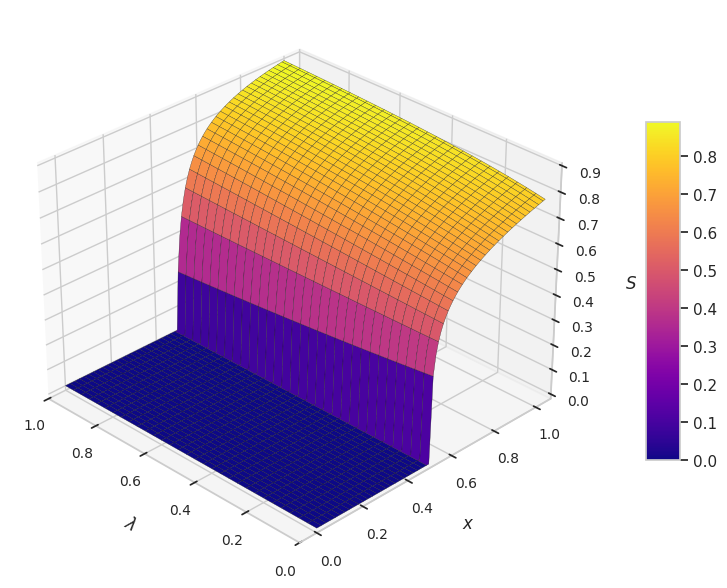}
            }
        \hfill
        \subcaptionbox{Convex/concave hull HC.}[0.40\linewidth]
            {
                \includegraphics[width=\linewidth]{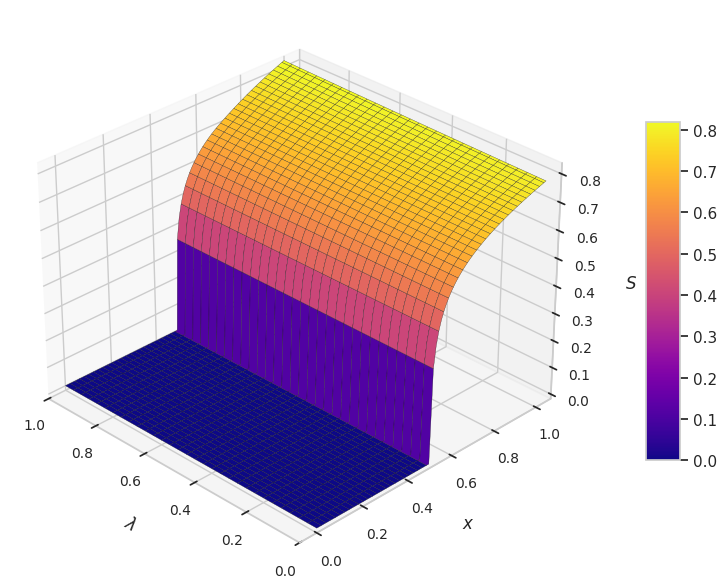}
            }
        \vspace{0cm}
        \subcaptionbox{Vanishing diffusion HC with \(\omega = 2 \times 10^{-3}\).}[0.40\linewidth]
            {
                \includegraphics[width=\linewidth]{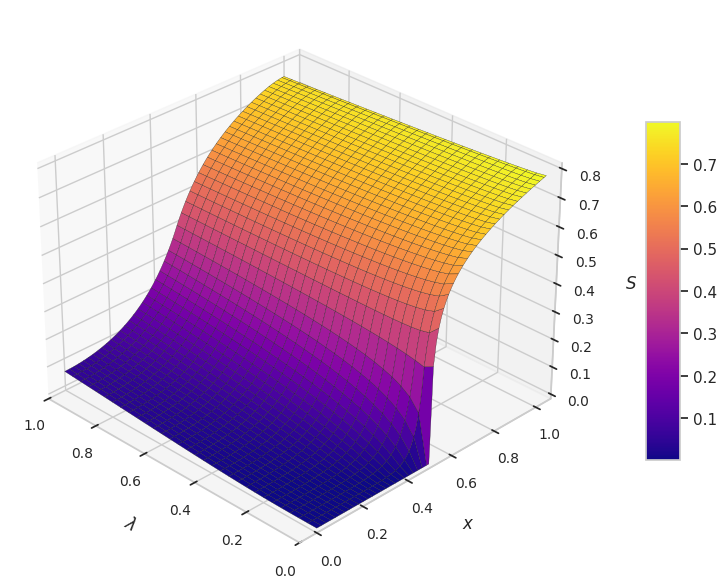}
            }
        \hfill
        \subcaptionbox{Vanishing diffusion HC with \(\omega = 10^{-5}\).}[0.40\linewidth]
            {
               \includegraphics[width=\linewidth]{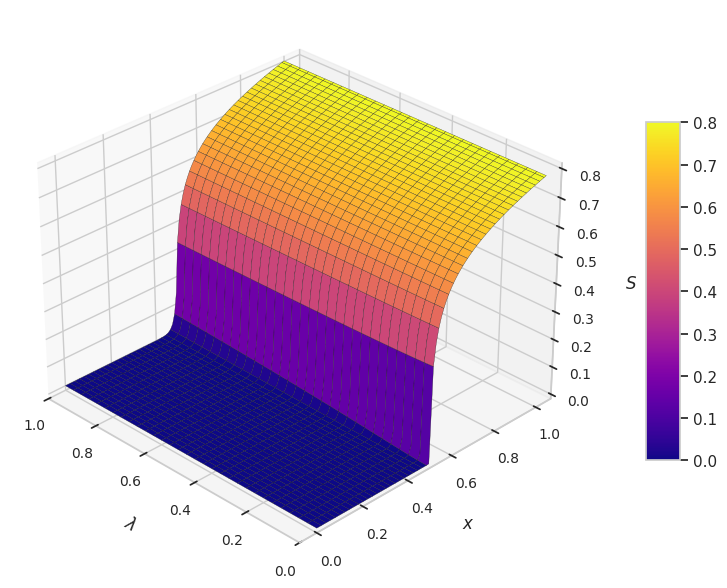}
            }
        \caption{Evolution of the solution along the HC curve beginning from the different auxiliary problems for case 2a.}
        \label{fig:results:case_2a:solution curves}
    \end{figure}

    The solution curves for Case 2a~\cref{fig:results:case_2a:solution curves} illustrate these trends. The convex/concave hull HC accurately reproduces the shock and the rarefaction wave of the target solution. The flow function in the linear relative permeability HC has a similar profile to the target flow function~\cref{fig:results:case_2b:flow function}, producing a similarly accurate solution despite lacking a shock. The vanishing diffusion HC with \(\omega = 2 \times 10^{-3}\) significantly smears the saturation front, whereas \(\omega = 10^{-5}\) yields an auxiliary solution that is almost indistinguishable from the target.

    Finally, cases symmetric with respect to phase labeling, for example the Riemann problem with \(M=0.1, S_\mr{left}^0 = 1.0, S_\mr{right}^0=0.0, S_i = 1.0\), which corresponds to case 2a, exhibit identical curvature and convergence behavior.

